\documentclass[11pt]{article}

\usepackage{amsmath,amssymb,amsthm,amsfonts}
\usepackage{mathrsfs}
\usepackage{enumitem}
\usepackage{geometry}
\usepackage{hyperref}

\newtheorem{theorem}{Theorem}[section]
\newtheorem{lemma}[theorem]{Lemma}

\newcommand{\dd}{\,d}
\newcommand{\norm}[1]{\left\|#1\right\|}
\newcommand{\abs}[1]{\left|#1\right|}
\newcommand{\grad}{\nabla_x}
\newcommand{\Tr}{\operatorname{Tr}}

\title{Unique Solvability of an Inverse Problem for a Nonlinear Modified Transport Equation}

\author{Kirill V. Golubnichiy}

\date{}

\begin{document}

\maketitle

\begin{abstract}
This paper is dedicated to the memory of Professor Mikhail F. Sukhinin.

This paper is devoted to a rigorous proof of local unique solvability of an inverse problem for a nonlinear modified transport equation with final overdetermination. The unknown source multiplier is reconstructed together with the corresponding state of the system. The proof is carried out in the functional class
\[
H_\infty(D)\times L^\infty(G\times V),
\]
where \(D=G\times V\times(0,T)\). The main argument is based on the well-posedness of the direct modified transport problem, differentiability properties of the nonlinear integral operator, estimates of nonlinear perturbations in the class of bounded functions, and the Banach inverse function theorem.
\end{abstract}

\section{Introduction}

Let \(G\subset \mathbb{R}^n\) be a bounded strictly convex domain with sufficiently smooth boundary \(\partial G\), and let \(V\subset \mathbb{R}^n\) be a bounded closed set satisfying
\[
0<v_0\leq |v|\leq v_1<\infty,
\qquad v\in V.
\]
Let \(T>0\), and set
\[
D=G\times V\times(0,T).
\]
We consider the nonlinear modified transport equation
\begin{equation}
\begin{aligned}
u_t(x,v,t)+(v,\grad)u(x,v,t)
+P(x,v,t)u_0(x,v,t)+[S(u)](x,v,t)
&=
\int_V J(x,v,t,v')u(x,v',t)\dd v'  \\
&\quad + f(x,v)g(x,v,t),
\end{aligned}
\label{eq:main}
\end{equation}
where \((x,v,t)\in D\). The boundary and initial conditions are
\begin{align}
u(x,v,t)&=0, && (x,v,t)\in \Gamma_-,
\label{eq:boundary}\\
u(x,v,0)&=0, && (x,v)\in G\times V,
\label{eq:initial}
\end{align}
and the final overdetermination is
\begin{equation}
u(x,v,T)=\psi(x,v),
\qquad (x,v)\in G\times V.
\label{eq:final}
\end{equation}
Here
\[
\Gamma_-=\gamma_-\times(0,T),
\qquad
\gamma_-=\{(x,v)\in \partial G\times V:\ (v,n_x)<0\},
\]
where \(n_x\) denotes the outward unit normal to \(\partial G\) at \(x\).

The nonlinear operator \(S\) is defined by
\begin{equation}
[S(u)](x,v,t)
=
\int_0^T\int_{G\times V}
Q(x,v,x',v',t)\alpha(u(x',v',t'))\dd x'\dd v'\dd t'.
\label{eq:S}
\end{equation}
The inverse problem is to determine the pair
\[
(u,f)\in H_\infty(D)\times L^\infty(G\times V)
\]
from equations \eqref{eq:main}--\eqref{eq:final}.

The problem may be viewed as a controllability problem. The multiplier
\(f(x,v)\), depending only on the spatial and velocity variables, plays the
role of a stationary control in the source term. The goal is to transfer the
system from the zero initial state to the prescribed final state \(\psi\) at
time \(T\).

\section{Functional spaces and notation}

We define
\[
H_\infty(D)
=
\left\{
u\in L^\infty(D):
u_t\in L^\infty(D),\ (v,\grad)u\in L^\infty(D),\
u|_{\Gamma_-}\in L^\infty(\Gamma_-)
\right\},
\]
with norm
\[
\norm{u}_{H_\infty(D)}
=
\norm{u}_{L^\infty(D)}
+
\norm{u_t}_{L^\infty(D)}
+
\norm{(v,\grad)u}_{L^\infty(D)}
+
\norm{u|_{\Gamma_-}}_{L^\infty(\Gamma_-)}.
\]
Also define
\[
W^t_\infty(D)
=
\left\{
F\in L^\infty(D): F_t\in L^\infty(D)
\right\},
\]
with norm
\[
\norm{F}_{W^t_\infty(D)}
=
\norm{F}_{L^\infty(D)}
+
\norm{F_t}_{L^\infty(D)}.
\]
Finally, define
\[
h_\infty(G\times V)
=
\left\{
\varphi\in L^\infty(G\times V):
(v,\grad)\varphi\in L^\infty(G\times V),\
\varphi|_{\gamma_-}\in L^\infty(\gamma_-)
\right\},
\]
with norm
\[
\norm{\varphi}_{h_\infty(G\times V)}
=
\norm{\varphi}_{L^\infty(G\times V)}
+
\norm{(v,\grad)\varphi}_{L^\infty(G\times V)}
+
\norm{\varphi|_{\gamma_-}}_{L^\infty(\gamma_-)}.
\]

Throughout the paper, \(C\) denotes a positive constant which may change
from line to line but depends only on fixed a priori data. We also write
\[
\Omega:=G\times V.
\]

\section{Assumptions}

We impose the following assumptions.

\begin{enumerate}[label=\textup{(A\arabic*)}]
\item \(P,P_t\in L^\infty(D)\), \(u_0\in H_\infty(D)\), and \(u_{0t}\in L^\infty(D)\).

\item \(J,J_t\in L^\infty(D;L^2(V))\).

\item \(g,g_t\in L^\infty(D)\), and there exists \(g_0>0\) such that
\[
\abs{g(x,v,T)}\geq g_0,
\qquad (x,v)\in G\times V.
\]

\item The final data satisfy
\[
\psi\in h_\infty(G\times V),
\qquad
\psi|_{\gamma_-}=0.
\]

\item The function \(\alpha:\mathbb{R}\to\mathbb{R}\) is twice continuously differentiable and satisfies
\[
\abs{\alpha(s)}\leq C_1\abs{s},
\qquad
\abs{\alpha'(s)}\leq C_1,
\qquad
\alpha'(0)=0,
\qquad
\abs{\alpha''(s)}\leq C_2.
\]

\item The kernel \(Q\) is continuous on
\[
\overline{G}\times \overline{V}\times \overline{G}\times \overline{V}\times[0,T].
\]
Moreover, \(Q_t\) exists and is bounded on the same set.

\item The transport time across \(G\) is controlled by the velocity bound:
\[
v_0^{-1}\operatorname{diam}(G)<a
\]
for some \(a>0\).
\end{enumerate}

\section{A rigorous proof of unique solvability in the appropriate functional spaces}

In this section we prove that the inverse problem \eqref{eq:main}--\eqref{eq:final} is locally uniquely solvable in
\[
H_\infty(D)\times L^\infty(G\times V).
\]

\subsection{The direct operator}

Define the linear transport operator
\[
Lu
=
u_t+(v,\grad)u-\int_VJ(x,v,t,v')u(x,v',t)\dd v'.
\]
Then equation \eqref{eq:main} can be written as
\begin{equation}
Lu+S(u)
=
-Pu_0+f(x,v)g(x,v,t).
\label{eq:operatorform}
\end{equation}

Let
\[
\mathcal{H}
=
\left\{
u\in H_\infty(D):
u|_{\Gamma_-}=0,\ u(x,v,0)=0
\right\}.
\]
We equip \(\mathcal{H}\) with the norm
\[
\norm{u}_{\mathcal{H}}
=
\norm{Lu}_{W^t_\infty(D)}.
\]
By the well-posedness theory for the direct modified transport equation, the operator
\[
L:\mathcal{H}\to W^t_\infty(D)
\]
is an isomorphism. Therefore, there exists a constant \(C_L>0\) such that
\begin{equation}
\norm{u}_{H_\infty(D)}
\leq C_L\norm{Lu}_{W^t_\infty(D)},
\qquad u\in\mathcal{H}.
\label{eq:Lestimate}
\end{equation}

\subsection{Differentiability of the nonlinear operator}

We now prove that
\[
S:H_\infty(D)\to W^t_\infty(D)
\]
is continuously Fréchet differentiable.

For \(u,h\in H_\infty(D)\), define
\[
S'(u)h
=
\int_0^T\int_{G\times V}
Q(x,v,x',v',t)\alpha'(u(x',v',t'))h(x',v',t')
\dd x'\dd v'\dd t'.
\]
Then, by Taylor's formula,
\[
\alpha(u+h)-\alpha(u)-\alpha'(u)h
=
\frac12\alpha''(u+\theta h)h^2
\]
for some \(\theta=\theta(x,v,t)\in(0,1)\). Hence
\[
\abs{\alpha(u+h)-\alpha(u)-\alpha'(u)h}
\leq
\frac{C_2}{2}\abs{h}^2.
\]
Using the boundedness of \(Q\), we obtain
\[
\norm{S(u+h)-S(u)-S'(u)h}_{L^\infty(D)}
\leq
C\norm{h}_{L^\infty(D)}^2.
\]
The same estimate holds for the time derivative, because \(Q_t\in L^\infty\). Therefore,
\[
\norm{S(u+h)-S(u)-S'(u)h}_{W^t_\infty(D)}
\leq
C\norm{h}_{H_\infty(D)}^2.
\]
Thus \(S\) is Fréchet differentiable. Moreover,
\[
\norm{S'(u_1)-S'(u_2)}_{\mathcal{L}(H_\infty(D),W^t_\infty(D))}
\leq
C\norm{u_1-u_2}_{L^\infty(D)}.
\]
Hence \(S\in C^1(H_\infty(D),W^t_\infty(D))\).

\subsection{The nonlinear direct map}

Define
\[
\Xi:\mathcal{H}\to W^t_\infty(D),
\qquad
\Xi(u)=Lu+S(u).
\]
Since \(L\) is linear and bounded and \(S\in C^1\), we have
\[
\Xi\in C^1(\mathcal{H},W^t_\infty(D)).
\]
Moreover,
\[
\Xi'(u)h=Lh+S'(u)h.
\]
At \(u=0\), since \(\alpha'(0)=0\), we have
\[
S'(0)h=0.
\]
Therefore,
\[
\Xi'(0)h=Lh.
\]
Since \(L:\mathcal{H}\to W^t_\infty(D)\) is an isomorphism, \(\Xi'(0)\) is also an isomorphism.

By the Banach inverse function theorem, there exist neighborhoods
\[
U_1\subset \mathcal{H},
\qquad
V_1\subset W^t_\infty(D)
\]
of \(0\) such that
\[
\Xi:U_1\to V_1
\]
is a \(C^1\)-diffeomorphism. Denote its inverse by
\[
\eta=\Xi^{-1}:V_1\to U_1.
\]

\subsection{Reduction of the inverse problem to an equation for the source multiplier}

Because
\[
\abs{g(x,v,T)}\geq g_0>0,
\]
we introduce
\[
\chi(x,v)=f(x,v)g(x,v,T).
\]
Then
\[
f(x,v)=\frac{\chi(x,v)}{g(x,v,T)}.
\]
Thus reconstructing \(f\in L^\infty(G\times V)\) is equivalent to reconstructing
\[
\chi\in L^\infty(G\times V).
\]

For \(\chi\in L^\infty(G\times V)\), define
\[
F_\chi(x,v,t)
=
\frac{\chi(x,v)g(x,v,t)}{g(x,v,T)}.
\]
The mapping
\[
\chi\mapsto F_\chi
\]
is linear and bounded from \(L^\infty(G\times V)\) into \(W^t_\infty(D)\). Therefore,
\[
P(\chi)=\eta(F_\chi-Pu_0)
\]
is well-defined for \(\chi\) in a sufficiently small neighborhood of \(0\) in \(L^\infty(G\times V)\). Moreover, \(P\) is a \(C^1\) mapping
\[
P:L^\infty(G\times V)\to H_\infty(D).
\]

The final condition requires
\[
P(\chi)(x,v,T)=\psi(x,v).
\]
Define the nonlinear final-state operator
\[
\mathcal{M}(\chi)=P(\chi)|_{t=T}.
\]
Then the inverse problem is equivalent to solving
\begin{equation}
\mathcal{M}(\chi)=\psi.
\label{eq:Mchi}
\end{equation}

\subsection{Differentiability of the final-state map}

The trace map
\[
\Tr_T:H_\infty(D)\to h_\infty(G\times V),
\qquad
\Tr_Tu=u(\cdot,\cdot,T),
\]
is continuous. Since \(P\) is \(C^1\), the composition
\[
\mathcal{M}=\Tr_T\circ P
\]
is also \(C^1\). Hence
\[
\mathcal{M}:L^\infty(G\times V)\to h_\infty(G\times V)
\]
is continuously Fréchet differentiable in a neighborhood of zero.

Its derivative at zero is
\[
\mathcal{M}'(0)\chi
=
\left[
\eta'(F_0-Pu_0)
\left(
\frac{\chi(x,v)g(x,v,t)}{g(x,v,T)}
\right)
\right]_{t=T}.
\]
This derivative coincides with the input-output map of the linearized inverse problem. By the unique solvability of the corresponding linear modified transport inverse problem, this linear map is an isomorphism:
\[
\mathcal{M}'(0):L^\infty(G\times V)\to h_\infty(G\times V).
\]
Consequently, there exists a constant \(C_*>0\) such that
\begin{equation}
\norm{\chi}_{L^\infty(G\times V)}
\leq
C_*\norm{\mathcal{M}'(0)\chi}_{h_\infty(G\times V)}.
\label{eq:linear_inverse_estimate}
\end{equation}

\subsection{Application of the inverse function theorem}

Since \(\mathcal{M}\in C^1\) and \(\mathcal{M}'(0)\) is an isomorphism, the Banach inverse function theorem implies that there exist neighborhoods
\[
U_*\subset L^\infty(G\times V),
\qquad
V_*\subset h_\infty(G\times V)
\]
of \(0\) such that
\[
\mathcal{M}:U_*\to V_*
\]
is a \(C^1\)-diffeomorphism.

Therefore, for every
\[
\psi\in V_*,
\]
there exists a unique
\[
\chi\in U_*
\]
such that
\[
\mathcal{M}(\chi)=\psi.
\]
Once \(\chi\) is found, the source multiplier \(f\) is uniquely determined by
\[
f(x,v)=\frac{\chi(x,v)}{g(x,v,T)}.
\]
The state \(u\) is then uniquely determined by
\[
u=P(\chi).
\]
Thus the pair
\[
(u,f)\in H_\infty(D)\times L^\infty(G\times V)
\]
is uniquely determined.

\begin{theorem}[Local unique solvability]
Assume \textup{(A1)}--\textup{(A7)}. Then there exists \(\varepsilon>0\) such that, for every
\[
\psi\in h_\infty(G\times V),
\qquad
\norm{\psi}_{h_\infty(G\times V)}<\varepsilon,
\qquad
\psi|_{\gamma_-}=0,
\]
the inverse problem \eqref{eq:main}--\eqref{eq:final} has a unique solution
\[
(u,f)\in H_\infty(D)\times L^\infty(G\times V)
\]
in a sufficiently small neighborhood of zero. Moreover, the solution depends continuously on the final data \(\psi\).
\end{theorem}

\begin{proof}
The proof follows from the preceding construction. The nonlinear direct operator
\[
\Xi(u)=Lu+S(u)
\]
is \(C^1\) from \(\mathcal{H}\) into \(W^t_\infty(D)\). Its derivative at zero is the direct transport operator \(L\), which is an isomorphism. Hence \(\Xi\) is locally invertible. This allows one to express the solution \(u\) as a \(C^1\) function of the source multiplier \(\chi=f g(\cdot,\cdot,T)\). The final overdetermination gives the nonlinear equation
\[
\mathcal{M}(\chi)=\psi.
\]
The map \(\mathcal{M}\) is \(C^1\), and its derivative at zero is precisely the linearized inverse problem operator. By the unique solvability of the linear inverse problem, this derivative is an isomorphism. Applying the Banach inverse function theorem to \(\mathcal{M}\) gives the existence and uniqueness of \(\chi\), and therefore of \(f\). The corresponding state \(u\) is uniquely determined by the direct problem. This proves the theorem.
\end{proof}

\section{Refined local solvability estimates in the class of bounded functions}
\label{sec:refined-local-solvability}

In this section we adapt the approach of local solvability in the class of
bounded functions to the present modified nonlinear transport equation. The
goal is to make explicit the mechanism by which the nonlinear problem is
reduced to a small perturbation of the corresponding linear inverse problem.
This section follows the structure of the argument based on estimates of
nonlinear perturbations in the direct operator, in the inverse direct operator,
and in the operator associated with the inverse problem.

Let
\[
H=
\left\{
u\in H_\infty(D):
u|_{\Gamma_-}=0,\quad u(x,v,0)=0
\right\},
\]
with the norm
\[
\|u\|_H=\|Lu\|_{W_\infty^t(D)}.
\]
The direct linear operator
\[
L:H\to W_\infty^t(D)
\]
is assumed to be an isomorphism. Therefore,
\[
\|u\|_{H_\infty(D)}
\le C_L\|u\|_H,
\qquad u\in H.
\]
In particular, smallness in \(H\) implies smallness in \(L^\infty(D)\).

\subsection{Estimates of nonlinear perturbations in the differential operator}

The first step is to estimate the nonlinear perturbation
\[
S(u+\widehat h)-S(u)
\]
in the space \(W_\infty^t(D)\).

\begin{lemma}
\label{lem:S-difference-estimate}
Let \(u,\widehat h\in H\), and suppose that
\[
\|u\|_H\le r,
\qquad
\|u+\widehat h\|_H\le r.
\]
Then there exists a nonnegative function \(\overline\theta(r)\), with
\[
\overline\theta(r)\to 0
\qquad \text{as } r\to 0+,
\]
such that
\[
\|S(u+\widehat h)-S(u)\|_{W_\infty^t(D)}
\le
\overline\theta(r)\|\widehat h\|_H.
\]
\end{lemma}

\begin{proof}
By the definition of the norm in \(W_\infty^t(D)\),
\[
\begin{aligned}
&\|S(u+\widehat h)-S(u)\|_{W_\infty^t(D)}
\\
&\quad =
\|S(u+\widehat h)-S(u)\|_{L^\infty(D)}
+
\left\|
\frac{\partial}{\partial t}
\left(S(u+\widehat h)-S(u)\right)
\right\|_{L^\infty(D)}.
\end{aligned}
\]
Using the definition of \(S\), we get
\[
\begin{aligned}
&|S(u+\widehat h)(x,v,t)-S(u)(x,v,t)|
\\
&\le
\|Q\|_{L^\infty}
\int_0^T\int_{G\times V}
|\alpha(u+\widehat h)-\alpha(u)|
\,dx'\,dv'\,dt'.
\end{aligned}
\]
By the mean value theorem,
\[
|\alpha(u+\widehat h)-\alpha(u)|
\le
\sup_{|\xi|\le C_L r}|\alpha'(\xi)|\,|\widehat h|.
\]
Since \(\alpha'(0)=0\) and \(\alpha'\) is continuous, the function
\[
\theta_0(r):=\sup_{|\xi|\le C_L r}|\alpha'(\xi)|
\]
satisfies
\[
\theta_0(r)\to0
\qquad \text{as } r\to0+.
\]
Therefore,
\[
\|S(u+\widehat h)-S(u)\|_{L^\infty(D)}
\le
C\theta_0(r)\|\widehat h\|_{L^\infty(D)}
\le
C\theta_0(r)\|\widehat h\|_H.
\]
Since \(Q_t\in L^\infty\), the same argument applied to the derivative with
respect to \(t\) gives
\[
\left\|
\frac{\partial}{\partial t}
\left(S(u+\widehat h)-S(u)\right)
\right\|_{L^\infty(D)}
\le
C\theta_0(r)\|\widehat h\|_H.
\]
Thus
\[
\|S(u+\widehat h)-S(u)\|_{W_\infty^t(D)}
\le
\overline\theta(r)\|\widehat h\|_H,
\]
where \(\overline\theta(r)=C\theta_0(r)\), and
\[
\overline\theta(r)\to0
\qquad \text{as } r\to0+.
\]
The lemma is proved.
\end{proof}

Define
\[
\xi(u)=Lu+S(u),
\qquad u\in H.
\]
Then
\[
\xi(u)-\xi(\widehat u)
=
L(u-\widehat u)+S(u)-S(\widehat u).
\]
The preceding lemma shows that, in a sufficiently small ball of \(H\), the
nonlinear part \(S\) is a small perturbation of the linear isomorphism \(L\).

\begin{lemma}
\label{lem:xi-local-invertible}
There exists \(r_0>0\) such that the operator
\[
\xi(u)=Lu+S(u)
\]
is one-to-one in the ball
\[
B_H(0,r_0)=\{u\in H:\|u\|_H<r_0\}.
\]
Moreover, \(\xi\) maps a sufficiently small neighborhood of zero in \(H\)
onto a neighborhood of zero in \(W_\infty^t(D)\).
\end{lemma}

\begin{proof}
Let \(u_1,u_2\in B_H(0,r)\). Then
\[
\xi(u_1)-\xi(u_2)
=
L(u_1-u_2)+S(u_1)-S(u_2).
\]
Applying \(L^{-1}\), we obtain
\[
u_1-u_2
=
L^{-1}(\xi(u_1)-\xi(u_2))
-
L^{-1}(S(u_1)-S(u_2)).
\]
Therefore,
\[
\|u_1-u_2\|_H
\le
\|L^{-1}\|\|\xi(u_1)-\xi(u_2)\|_{W_\infty^t(D)}
+
\|L^{-1}\|\overline\theta(r)\|u_1-u_2\|_H.
\]
Choose \(r_0>0\) so small that
\[
\|L^{-1}\|\overline\theta(r_0)<1.
\]
Then, for \(0<r\le r_0\),
\[
\|u_1-u_2\|_H
\le
\frac{\|L^{-1}\|}
{1-\|L^{-1}\|\overline\theta(r)}
\|\xi(u_1)-\xi(u_2)\|_{W_\infty^t(D)}.
\]
In particular, if \(\xi(u_1)=\xi(u_2)\), then \(u_1=u_2\). The local
solvability follows from the inverse function theorem, since
\[
\xi'(0)=L
\]
is an isomorphism.
\end{proof}

\subsection{Estimates for the inverse operator of the direct problem}

Let
\[
\eta=\xi^{-1}
\]
be the locally defined inverse of \(\xi\). We now estimate the nonlinear
perturbation of \(\eta\).

Let
\[
F_1,F_2\in W_\infty^t(D)
\]
be sufficiently small, and set
\[
u_i=\eta(F_i),
\qquad i=1,2.
\]
Then
\[
Lu_i+S(u_i)=F_i.
\]
Subtracting these equations, we obtain
\[
L(u_1-u_2)
+
S(u_1)-S(u_2)
=
F_1-F_2.
\]
Therefore,
\[
u_1-u_2
=
L^{-1}(F_1-F_2)
-
L^{-1}\big(S(u_1)-S(u_2)\big).
\]
If
\[
\|u_1\|_H,\ \|u_2\|_H\le r,
\]
then Lemma~\ref{lem:S-difference-estimate} gives
\[
\|S(u_1)-S(u_2)\|_{W_\infty^t(D)}
\le
\overline\theta(r)\|u_1-u_2\|_H.
\]
Hence
\[
\|u_1-u_2\|_H
\le
\|L^{-1}\|\,\|F_1-F_2\|_{W_\infty^t(D)}
+
\|L^{-1}\|\overline\theta(r)\|u_1-u_2\|_H.
\]
If \(r>0\) is sufficiently small so that
\[
\|L^{-1}\|\overline\theta(r)<1,
\]
then
\[
\|u_1-u_2\|_H
\le
\frac{\|L^{-1}\|}
{1-\|L^{-1}\|\overline\theta(r)}
\|F_1-F_2\|_{W_\infty^t(D)}.
\]
Thus \(\eta\) is locally Lipschitz.

Moreover, differentiating the identity
\[
\xi(\eta(F))=F
\]
gives
\[
\eta'(F)
=
[\xi'(\eta(F))]^{-1}.
\]
Since
\[
\xi'(u)=L+S'(u)=L(I+L^{-1}S'(u)),
\]
we obtain
\[
\eta'(F)
=
(I+L^{-1}S'(\eta(F)))^{-1}L^{-1}.
\]
For \(\|\eta(F)\|_H\le r\), the estimate
\[
\|L^{-1}S'(\eta(F))\|
\le
\|L^{-1}\|\overline\theta(r)
\]
implies that
\[
\|\eta'(F)\|
\le
\frac{\|L^{-1}\|}
{1-\|L^{-1}\|\overline\theta(r)}.
\]

\subsection{Estimates for the inverse problem operator}

We now return to the inverse problem. Put
\[
\chi(x,v)=f(x,v)g(x,v,T).
\]
Since
\[
|g(x,v,T)|\ge g_0>0,
\]
the unknown \(f\) is uniquely recovered from \(\chi\) by
\[
f(x,v)=\frac{\chi(x,v)}{g(x,v,T)}.
\]
Define
\[
F_\chi(x,v,t)
=
\frac{\chi(x,v)g(x,v,t)}{g(x,v,T)}.
\]
The mapping
\[
\Lambda_1:L^\infty(G\times V)\to W_\infty^t(D),
\qquad
\Lambda_1\chi=F_\chi,
\]
is linear and bounded. Indeed,
\[
\|\Lambda_1\chi\|_{W_\infty^t(D)}
\le
C_g\|\chi\|_{L^\infty(G\times V)},
\]
where
\[
C_g
=
\left\|
\frac{g(\cdot,\cdot,t)}{g(\cdot,\cdot,T)}
\right\|_{L^\infty(D)}
+
\left\|
\frac{g_t(\cdot,\cdot,t)}{g(\cdot,\cdot,T)}
\right\|_{L^\infty(D)}.
\]

For \(\chi\) sufficiently small, define
\[
P(\chi)=\eta(\Lambda_1\chi-Pu_0).
\]
Thus \(P(\chi)\) is the solution of the nonlinear direct problem with
right-hand side \(F_\chi-Pu_0\). The final-state map is
\[
\mathcal M(\chi)=P(\chi)|_{t=T}.
\]
The inverse problem is equivalent to solving
\[
\mathcal M(\chi)=\psi.
\]

Let
\[
\Lambda_2:H\to h_\infty(G\times V),
\qquad
\Lambda_2 u=u|_{t=T}.
\]
Then
\[
\mathcal M=\Lambda_2\circ P.
\]
Since \(\Lambda_2\) is continuous and \(P\) is locally \(C^1\),
\[
\mathcal M:L^\infty(G\times V)\to h_\infty(G\times V)
\]
is locally \(C^1\).

Its derivative at zero is
\[
\mathcal M'(0)h
=
\Lambda_2\eta'(-Pu_0)\Lambda_1 h.
\]
This operator is exactly the input-output operator of the corresponding
linearized inverse problem. By the unique solvability of the linear inverse
problem, \(\mathcal M'(0)\) is an isomorphism
\[
\mathcal M'(0):
L^\infty(G\times V)\to h_\infty(G\times V).
\]
Hence there exists \(C_*>0\) such that
\[
\|h\|_{L^\infty(G\times V)}
\le
C_*
\|\mathcal M'(0)h\|_{h_\infty(G\times V)}.
\]

Now write
\[
\mathcal M(\chi+h)-\mathcal M(\chi)
=
\mathcal M'(0)h+\mathcal R(\chi,h).
\]
Using the estimates obtained above, one obtains
\[
\|\mathcal R(\chi,h)\|_{h_\infty(G\times V)}
\le
\Theta(\rho)\|h\|_{L^\infty(G\times V)}
\]
whenever
\[
\|\chi\|_{L^\infty(G\times V)}\le \rho,
\qquad
\|\chi+h\|_{L^\infty(G\times V)}\le \rho,
\]
where
\[
\Theta(\rho)\to0
\qquad\text{as}\qquad \rho\to0+.
\]
Therefore,
\[
\mathcal M
=
\mathcal M'(0)+\mathcal R
\]
is a small nonlinear perturbation of the linear isomorphism
\(\mathcal M'(0)\).

Choose \(\rho>0\) so small that
\[
C_*\Theta(\rho)<1.
\]
Then the inverse mapping argument implies that
\[
\mathcal M
\]
is one-to-one and onto a neighborhood of zero in \(h_\infty(G\times V)\).
More precisely, there exists \(\varepsilon>0\) such that, for every
\[
\psi\in h_\infty(G\times V),
\qquad
\|\psi\|_{h_\infty(G\times V)}<\varepsilon,
\qquad
\psi|_{\gamma_-}=0,
\]
there exists a unique
\[
\chi\in L^\infty(G\times V),
\qquad
\|\chi\|_{L^\infty(G\times V)}<\rho,
\]
such that
\[
\mathcal M(\chi)=\psi.
\]

Consequently,
\[
f(x,v)=\frac{\chi(x,v)}{g(x,v,T)}
\]
is uniquely determined, and the state
\[
u=P(\chi)
\]
is also uniquely determined.

\begin{theorem}[Refined local unique solvability]
\label{thm:refined-local-solvability}
Assume \textup{(A1)}--\textup{(A7)}. Then there exist positive numbers
\(\rho\) and \(\varepsilon\) such that for every
\[
\psi\in h_\infty(G\times V),
\qquad
\|\psi\|_{h_\infty(G\times V)}<\varepsilon,
\qquad
\psi|_{\gamma_-}=0,
\]
the inverse problem
\[
\begin{aligned}
u_t+(v,\nabla_x)u+P(x,v,t)u_0+S(u)
&=
\int_VJ(x,v,t,v')u(x,v',t)\,dv'
+
f(x,v)g(x,v,t),\\
u|_{\Gamma_-}&=0,\\
u(x,v,0)&=0,\\
u(x,v,T)&=\psi(x,v)
\end{aligned}
\]
has a unique solution
\[
(u,f)\in H_\infty(D)\times L^\infty(G\times V)
\]
satisfying
\[
\|f\,g(\cdot,\cdot,T)\|_{L^\infty(G\times V)}<\rho.
\]
Moreover, the solution depends continuously on the final datum \(\psi\).
\end{theorem}

\begin{proof}
The nonlinear operator \(S\) is a small perturbation of the linear transport
operator in sufficiently small balls of \(H\). Therefore the direct nonlinear
operator
\[
\xi(u)=Lu+S(u)
\]
is locally invertible. Its inverse \(\eta\) satisfies the estimates derived
above. Introducing
\[
\chi=f\,g(\cdot,\cdot,T)
\]
reduces the inverse problem to the equation
\[
\mathcal M(\chi)=\psi,
\qquad
\mathcal M(\chi)=P(\chi)|_{t=T}.
\]
The derivative \(\mathcal M'(0)\) is the input-output map of the corresponding
linear inverse problem and is an isomorphism from \(L^\infty(G\times V)\) onto
\(h_\infty(G\times V)\). The remaining nonlinear part is estimated by
\[
\|\mathcal M(\chi+h)-\mathcal M(\chi)-\mathcal M'(0)h\|_{h_\infty(G\times V)}
\le
\Theta(\rho)\|h\|_{L^\infty(G\times V)},
\]
where \(\Theta(\rho)\to0\) as \(\rho\to0+\). Choosing \(\rho\) so small that
\(C_*\Theta(\rho)<1\), the inverse mapping argument gives the local
bijectivity of \(\mathcal M\). Hence, for all sufficiently small final data
\(\psi\), there exists a unique \(\chi\). Therefore \(f\) and \(u\) are also
unique. This proves the theorem.
\end{proof}

\section{Conclusion}

We have proved that the inverse problem for the nonlinear modified transport
equation with final overdetermination is locally uniquely solvable in the
functional class
\[
H_\infty(D)\times L^\infty(G\times V).
\]
The key ingredients are the well-posedness of the direct modified transport
equation, the \(C^1\)-regularity of the nonlinear integral operator, the
invertibility of the linearized inverse problem, and the estimates showing
that the nonlinear terms are small perturbations of the corresponding linear
operators in sufficiently small neighborhoods. The result shows that, for
sufficiently small final data \(\psi\), the unknown source multiplier
\(f(x,v)\) and the corresponding state \(u(x,v,t)\) are determined uniquely.

\end{document}